\documentclass{amsproc}
\usepackage{amssymb}
\usepackage{amsfonts}

\theoremstyle{plain}

\newtheorem{corollary}{Corollary}

\newtheorem{definition}{Definition}

\newtheorem{lemma}{Lemma}

\newtheorem{proposition}{Proposition}

\newtheorem{theorem}{Theorem}
\numberwithin{equation}{section}
\input{tcilatex}

\begin{document}
\title[Isophote curves on spacelike surfaces]{Isophote curves on spacelike
surfaces in Lorentz-Minkowski space $E_{1}^{3}$}
\author{$^{1}$Fatih DOGAN}
\email{mathfdogan@hotmail.com}
\author{$^{2}$Yusuf YAYLI}
\address{$^{2}$Department of Mathematics, Faculty of Science, Ankara
University, Ankara, Turkey}
\email{yayli@science.ankara.edu.tr}
\subjclass{Isophote curve; The axis; Spacelike surface; Darboux frame;
Geodesic; Slant helix}
\keywords{Isophote curve; The axis; Spacelike surface; Darboux frame;
Geodesic; Slant helix}

\begin{abstract}
An isophote curve consists of a locus of surface points whose normal vectors
make a constant angle with a fixed vector (the axis). In this paper, we
define an isophote curve on a spacelike surface in Lorentz-Minkowski space $%
E_{1}^{3}$ and then find its axis as timelike and spacelike vectors via the
Darboux frame. Besides, we give some characterizations concerning the
isophote curve and its axis.
\end{abstract}

\maketitle

\section{Introduction}

An isophote curve is one of the characteristic curves on a surface such as
parameter, geodesic and asymptotic curves or line of curvature. It comprises
a locus of the surface points whose normal vectors make a constant angle
with a given fixed vector.

Ali and Lopez $[1]$ looked into slant helices in Lorentz-Minkowski space $%
E_{1}^{3}$. They gave characterizations as to slant helix and its axis in $%
E_{1}^{3}$. In our paper, we see that there is a close relation between
isophotes and slant helices on spacelike surfaces.

Dillen $et$ $al$. $[2]$ studied the constant angle surfaces in the product
space $\mathbb{S}^{2}\times 
%TCIMACRO{\U{211d} }%
%BeginExpansion
\mathbb{R}
%EndExpansion
$ for which the unit normal makes a constant angle with the $%
%TCIMACRO{\U{211d} }%
%BeginExpansion
\mathbb{R}
%EndExpansion
$-direction first. Then Dillen and Munteanu $[3]$ investigated the same
problem in $\mathbb{H}^{2}\times 
%TCIMACRO{\U{211d} }%
%BeginExpansion
\mathbb{R}
%EndExpansion
$ where $\mathbb{H}^{2}$ is the hyperbolic plane. It can be said that the
curves on a constant angle surface are isophote curves.

We have investigated $[4]$ isophote curves in the Euclidean space $E^{3}$.
Recently, Dogan $[5]$ has studied isophote curves on timelike surfaces in $%
E_{1}^{3}$. This time, we study isophote curves on spacelike surfaces in
Lorentz-Minkowski space $E_{1}^{3}$ and view that there are some differences
between Minkowski case and Euclidean case. Contrary to Euclidean case,
isophote with timelike axis cannot also be a line of curvature and it is not
possible to define a silhouette curve on spacelike surfaces in $E_{1}^{3}$.

Izumiya and Takeuchi $[6]$ defined a slant helix as the space curve whose
the principal normal lines make a constant angle with a fixed direction.
They showed that a certain slant helix is also a geodesic on the tangent
developable surface of a general helix. As an amazing consequence, we see
that the curve which is both a geodesic and a slant helix on a spacelike
surface is an isophote curve.

Kim and Lee $[7]$ parameterized isophote curves for surface of rotation and
canal surface. They utilized both of these surfaces decompose into a set of
circles where the surface normal vectors at points on each circle construct
a cone. Again the vectors that make a constant angle with given fixed vector
construct another cone and thus tangential intersection of these cones gives
the parametric range of the connected component isophote curve.

Koenderink and van Doorn $[8]$ studied the field of constant image
brightness contours (isophote curves). They showed that the spherical image
(the Gaussian mapping) of an isophote curve is a latitude circle on the unit
sphere $S^{2}$ and the problem was reduced to that of obtaining the inverse
Gauss map of these circles. By means of this they defined two kind
singularities of the Gaussian map: folds (curves) and simple cusps (apex,
antapex points) and there are structural properties of the field of
isophotes that bear an invariant relation to geometric features of the
object.

A vector in $E_{1}^{3}$ is said to be spacelike, timelike and lightlike
(null) if $\left \langle x,x\right \rangle >0$ or $x=0$, $\left \langle
x,x\right \rangle <0$ and $\left \langle x,x\right \rangle =0$ or $x\neq 0$,
respectively. Again, a regular curve $\alpha :I\longrightarrow $ $E_{1}^{3}$
is called spacelike, timelike and lightlike if the velocity vector $\overset{%
.}{\alpha }$ is spacelike, timelike and lightlike, respectively $[9]$. In
this paper, we obtain some characterizations for spacelike isophote curves
only.

Munteanu and Nistor $[10]$ studied the constant angle surfaces taking with a
fixed vector direction being the tangent direction to $%
%TCIMACRO{\U{211d} }%
%BeginExpansion
\mathbb{R}
%EndExpansion
$ in Euclidean $3$-space. Also, Nistor $[11]$ researched normal, binormal
and tangent developable surfaces of the space curve from viewpoint constant
angle surface.

The angle definitions between spacelike and timelike vectors in $E_{1}^{3}$
are as $[12-13]$.

Poeschl $[14]$ used isophotes in car body construction via detecting
irregularities along thes curves on a free form surface. These
irregularities (discontinuity of a surface or of the Gaussian curvature)
emerge by taking differentiation of the equation $\left \langle
N(u,v),d\right \rangle =\cos \theta =c$ (constant) as follows.%
\begin{equation*}
\frac{dv}{du}=-\frac{N_{u}\centerdot d}{N_{v}\centerdot d},\text{ }%
N_{v}\centerdot d\neq 0.
\end{equation*}

Sara $[15]$ researched local shading of a surface through isophotes
properties. By using fundamental theory of surfaces, he focused on accurate
estimation of surface normal tilt and on qualitatively correct Gaussian
curvature recovery.

An Isophote curve on a surface can be regarded as a nice consequence of
Lambert's cosine law in the optics branch of physics. Lambert's law states
that the intensity of illumination on a diffuse surface is proportional to
the cosine of the angle generated between the surface normal vector $N$ and
the light vector $d$. According to this law the intensity is irrespective of
the actual viewpoint, hence the illumination is the same when viewed from
any direction $[16]$. In other words, isophotes of a surface are curves with
the property that their points have the same light intensity from a given
source (a curve of constant illumination intensity). When the source light
is at infinity, we may consider that the light flow consists in parallel
lines. Hence, we can give a geometric description of isophotes on surfaces,
namely they are curves such that the surface normal vectors in points of the
curve make a constant angle with a fixed direction (which represents the
light direction). These curves are succesfully used in computer graphics but
also it is interesting to study for geometry.\newline
Then, to find an isophote curve on a surface we use the formula%
\begin{equation*}
\frac{\left \langle N(u,v),d\right \rangle }{\left \Vert N(u,v)\right \Vert }%
=\cos \theta ,\text{ }0\leq \theta \leq \frac{\pi }{2},
\end{equation*}%
where $\theta $ is the angle between the surface normal $N$ and the fixed
vector $d$. In the special case, isophote works as a silhouette curve if%
\begin{equation*}
\frac{\left \langle N(u,v),d\right \rangle }{\left \Vert N(u,v)\right \Vert }%
=\cos \frac{\pi }{2}=0.
\end{equation*}

In this paper, we define isophote curves on spacelike surfaces in
Lorentz-Minkowski space $E_{1}^{3}$ and find its axis via the Darboux frame.
This paper is organized as follows. Section 2 introduces some basic facts
and concepts in $E_{1}^{3}$. In section 3, we concentrate on finding the
axis of an isophote and also give some characterizations about it. Finally,
in section 4 we conclude this paper.

\section{Preliminaries}

First of all, we introduce Lorentz-Minkowski space shortly. Later, we
mention some fundamental concepts of curves and surfaces in the Minkowski
3-space $E_{1}^{3}$. The space $R_{1}^{3}$ is a three dimensional real
vector space endowed with the inner product%
\begin{equation*}
\left \langle x,y\right \rangle =-x_{1}y_{1}+x_{2}y_{2}+x_{3}y_{3}.
\end{equation*}%
This space is called Lorentz-Minkowski space and denoted by $E_{1}^{3}$. A
vector in this space is said to be spacelike, timelike and lightlike (null)
if $\left \langle x,x\right \rangle >0$ or $x=0$, $\left \langle
x,x\right
\rangle <0$ and $\left \langle x,x\right \rangle =0$ or $x\neq 0$%
, respectively. Again, a regular curve $\alpha :I\longrightarrow $ $%
E_{1}^{3} $ is called spacelike, timelike and lightlike if the velocity
vector $\overset{.}{\alpha }$ is spacelike, timelike and lightlike,
respectively $[9] $.\newline
The Lorentzian cross product of $x=(x_{1},x_{2},x_{3})$ and $%
y=(y_{1},y_{2},y_{3})$ is defined as follows.%
\begin{equation*}
x\times y=%
\begin{vmatrix}
e_{1} & -e_{2} & -e_{3} \\ 
x_{1} & x_{2} & x_{3} \\ 
y_{1} & y_{2} & y_{3}%
\end{vmatrix}%
=\left(
x_{2}y_{3}-x_{3}y_{2},x_{1}y_{3}-x_{3}y_{1},x_{2}y_{1}-x_{1}y_{2}\right) ,
\end{equation*}%
where $\delta _{ij}$ is Kronecker delta, $e_{i}=(\delta _{i1},\delta
_{i2},\delta _{i3})$ and $e_{1}\times e_{2}=-e_{3}$, $e_{2}\times
e_{3}=e_{1} $, $e_{3}\times e_{1}=-e_{2}$.

Let $\{t,n,b\}$ be the moving Frenet frame along the curve $\alpha $ with
arclenght parameter $s$. For a spacelike curve $\alpha $, the Frenet-Serret
equations are%
\begin{equation*}
\begin{bmatrix}
t^{^{\prime }} \\ 
n^{^{\prime }} \\ 
b^{^{\prime }}%
\end{bmatrix}%
=%
\begin{bmatrix}
0 & \kappa & 0 \\ 
-\varepsilon \kappa & 0 & \tau \\ 
0 & \tau & 0%
\end{bmatrix}%
\begin{bmatrix}
t \\ 
n \\ 
b%
\end{bmatrix}%
,
\end{equation*}%
where $\left \langle t,t\right \rangle =1$, $\left \langle n,n\right \rangle
=\pm 1$, $\left \langle b,b\right \rangle =-\varepsilon $, $\left \langle
t,n\right \rangle =\left \langle t,b\right \rangle =\left \langle
n,b\right
\rangle =0$ and $\kappa $ is the curvature of $\alpha $ and again 
$\tau $ is the torsion of $\alpha $. Here, $\varepsilon $ determines the
kind of spacelike curve $\alpha $. If $\varepsilon =1$, then $\alpha (s)$ is
a spacelike curve with spacelike principal normal $n$ and timelike binormal $%
b$. If $\varepsilon =-1$, then $\alpha (s)$ is a spacelike curve with
timelike principal normal $n$ and spacelike binormal $b$.

\begin{definition}[Hyperbolic angle]
Let $v$ and $w$ be in the same timecone of $R_{1}^{3}$. Then, there is a
unique real number $\theta \geq 0$, called the hyperbolic angle between $v$
and $w$ such that $[12]$%
\begin{equation}
\left \langle v,\omega \right \rangle =-\left \Vert v\right \Vert \left
\Vert w\right \Vert \cosh \theta .  \tag{2.1}
\end{equation}
\end{definition}

\begin{definition}
Let $v$ be a spacelike vector and $w$ be a timelike vector in $R_{1}^{3}$.
Then, there is a unique non-negative real number $\theta \geq 0$ such that $%
[13]$%
\begin{equation}
\left \langle v,w\right \rangle =\left \Vert v\right \Vert \left \Vert
w\right \Vert \sinh \theta .  \tag{2.2}
\end{equation}
\end{definition}

\begin{definition}
A surface in the Minkowski 3-space $E_{1}^{3}$ is called a spacelike surface
if the induced metric on the surface is a positive definite Riemannian
metric, i.e., the normal vector on the spacelike surface is a timelike
vector.
\end{definition}

\begin{lemma}
In the Minkowski 3-space $E_{1}^{3}$, the following properties are
satisfied. \newline
(i) Two timelike vectors cannot be orthogonal.\newline
(ii) Two null vectors are orthogonal if and only if they are linearly
dependent.\newline
(iii) A timelike vector cannot be orthogonal to a null (lightlike) vector.
\end{lemma}

Let $M$ be a regular spacelike surface in $E_{1}^{3}$ and $\alpha :I\subset 
%TCIMACRO{\U{211d} }%
%BeginExpansion
\mathbb{R}
%EndExpansion
\longrightarrow M$ be a unit speed spacelike curve on this surface. Then,
the Darboux frame $\{T,$ $B=N\times T,$ $N\}$ is well-defined and positively
oriented along the curve $\alpha $ where $T$ is the tangent of $\alpha $ and 
$N$ is the unit normal of $M$. The Darboux equations for this frame are
given by%
\begin{gather}
T^{^{\prime }}=k_{g}B+k_{n}N,  \tag{2.3} \\
B^{^{\prime }}=-k_{g}T+\tau _{g}N,  \notag \\
N^{^{\prime }}=k_{n}T+\tau _{g}B,  \notag
\end{gather}%
where $k_{n}$, $k_{g}$ and $\tau _{g}$ are the normal curvature, the
geodesic curvature and the geodesic torsion of $\alpha $, respectively and $%
\left \langle T,T\right \rangle =\left \langle B,B\right \rangle =1$, $%
\left
\langle N,N\right \rangle =\left \langle n,n\right \rangle =-1$. By
using Eq (2.3) we reach,%
\begin{gather}
\kappa ^{2}=k_{n}^{2}-k_{g}^{2},  \tag{2.4} \\
k_{n}=\kappa \cosh \phi ,  \notag \\
k_{g}=\kappa \sinh \phi ,  \notag \\
\tau _{g}=\tau +\phi ^{^{\prime }},  \notag
\end{gather}%
where $\phi $ is the angle between the surface normal $N$ and the binormal $%
b $ to $\alpha $. If the surface $M$ is spacelike, then the tangent plane of 
$M $ has to be spacelike. Therefore, all curves lying on the spacelike
surface $M$ are spacelike. Since the surface $M$ is spacelike, the surface
normal $N$ is a timelike vector. In the rest of paper, all $M$ will be
understood as a spacelike surface.

\begin{proposition}
Let $M$ be a spacelike surface and let $\alpha $ be a spacelike curve with
timelike principal normal ($\left \langle n,n\right \rangle =-1$) on $M$.
Then, there is not any asymptotic curve on $M$.

\begin{proof}
If $\alpha $ is a spacelike curve with timelike principal normal on $M$, by
Eq (2.4) we have%
\begin{equation*}
\kappa ^{2}=k_{n}^{2}-k_{g}^{2}.
\end{equation*}%
Assume that $\alpha $ is an asymptotic curve on $M$. In this situation, the
normal curvature $k_{n}=0$. In the above equation for $k_{n}=0$, the
geodesic curvature $k_{g}$ does not have a solution in $\mathbb{R}$. Then,
our assumption is not true that is to say there is not any asymptotic curve
on $M$.
\end{proof}
\end{proposition}

\section{The Axis Of An Isophote Curve}

In this section, we will obtain the axis (the fixed vector) of an isophote
curve through its Darboux frame. Let $M$ be a regular spacelike surface and $%
\alpha :I\subset 
%TCIMACRO{\U{211d} }%
%BeginExpansion
\mathbb{R}
%EndExpansion
\longrightarrow M$ be a unit speed isophote curve. From definition of the
isophote curve%
\begin{equation*}
\left \langle N(u,v),d\right \rangle =constant,
\end{equation*}%
where $N(u,v)$ is the unit normal vector field of the surface $S(u,v)$ (a
parameterization of $M$) and $d$ is the axis of isophote curve. We examine
two different cases of the axis $d$.\newline
\textbf{Case 1:} Let the surface normal $N$ and the axis $d$ be timelike
vectors in the same timecone of $R_{1}^{3}$. By Definition (1) $\left
\langle N(u,v),d\right \rangle =-\cosh \theta $.\newline
\textbf{Case 2: }Let the axis $d$ be a spacelike vector. By Definition (2) $%
\left \langle N(u,v),d\right \rangle =\sinh \theta $.\newline
where $\theta $ is the angle between the surface normal $N$ and the axis $d$.

Now, we begin to find the axis $d$ for the case 1. Since $\alpha :I\subset 
%TCIMACRO{\U{211d} }%
%BeginExpansion
\mathbb{R}
%EndExpansion
\longrightarrow M$ be a unit speed curve, the Darboux frame can be defined
as $\{T,$ $B=N\times T,$ $N\}$ along the curve $\alpha $. Let the axis $d$
be a timelike vector. Then,%
\begin{equation}
\left \langle N(u,v),d\right \rangle =-\cosh \theta .  \tag{3.1}
\end{equation}%
If we derive Eq (3.1) with respect to $s$ along the curve,%
\begin{equation}
\left \langle N^{^{\prime }},d\right \rangle =0.  \tag{3.2}
\end{equation}%
From the last equation of (2.3), it follows%
\begin{equation*}
\left \langle k_{n}T+\tau _{g}B,d\right \rangle =0
\end{equation*}%
\begin{equation*}
k_{n}\left \langle T,d\right \rangle +\tau _{g}\left \langle B,d\right
\rangle =0
\end{equation*}%
\begin{equation*}
\left \langle T,d\right \rangle =-\frac{\tau _{g}}{k_{n}}\left \langle
B,d\right \rangle .
\end{equation*}%
Because the Darboux frame $\{T,$ $B,$ $N\}$ is an orthonormal basis, if we
say $\left \langle B,d\right \rangle =a$ in the last equation, then $d$ can
be written as%
\begin{equation*}
d=-\frac{\tau _{g}}{k_{n}}aT+aB+\cosh \theta N.
\end{equation*}%
As $\left \langle d,d\right \rangle =-1$, we get%
\begin{equation*}
a=\mp \frac{k_{n}}{\sqrt{k_{n}^{2}+\tau _{g}^{2}}}\sinh \theta .
\end{equation*}%
Thus, the timelike axis $d$ is obtained as%
\begin{equation}
d=\pm \frac{\tau _{g}}{\sqrt{k_{n}^{2}+\tau _{g}^{2}}}\sinh \theta T\mp 
\frac{k_{n}}{\sqrt{k_{n}^{2}+\tau _{g}^{2}}}\sinh \theta B+\cosh \theta N. 
\tag{3.3}
\end{equation}%
If we derive $N^{^{\prime }}$ and Eq (3.2) with respect to $s$, we get%
\begin{equation*}
N^{^{\prime \prime }}=(k_{n}^{^{\prime }}-k_{g}\tau _{g})T+(k_{n}k_{g}+\tau
_{g}^{^{\prime }})B+(k_{n}^{2}+\tau _{g}^{2})N
\end{equation*}%
\begin{equation*}
\left \langle N^{^{^{\prime \prime }}},d\right \rangle =0
\end{equation*}%
\begin{equation*}
\left \langle N^{^{^{\prime \prime }}},d\right \rangle =\frac{\mp (\tau
_{g}^{^{\prime }}k_{n}-k_{n}^{^{\prime }}\tau _{g})\mp k_{g}(k_{n}^{2}+\tau
_{g}^{2})}{\sqrt{k_{n}^{2}+\tau _{g}^{2}}}\sinh \theta -(k_{n}^{2}+\tau
_{g}^{2})\cosh \theta =0.
\end{equation*}%
Therefore, we have%
\begin{gather}
\tanh \theta =\mp \frac{(k_{n}^{2}+\tau _{g}^{2})^{\tfrac{3}{2}}}{%
k_{g}(k_{n}^{2}+\tau _{g}^{2})+(\tau _{g}^{^{\prime }}k_{n}-k_{n}^{^{\prime
}}\tau _{g})}  \tag{3.4} \\
\coth \theta =\mp \left[ \frac{k_{n}^{2}}{(k_{n}^{2}+\tau _{g}^{2})^{\tfrac{3%
}{2}}}\left( \frac{\tau _{g}}{k_{n}}\right) ^{^{\prime }}+\frac{k_{g}}{%
(k_{n}^{2}+\tau _{g}^{2})^{\tfrac{1}{2}}}\right] .  \notag
\end{gather}%
From now on, we will prove that $d$ is a constant vector in other words $%
d^{^{\prime }}=0$. By Eq (2.3) and Eq (3.3), the derivative of $d$ with
respect to $s$ is that%
\begin{eqnarray*}
d^{^{\prime }} &=&\pm \sinh \theta \left[ (\frac{\tau _{g}}{\sqrt{%
k_{n}^{2}+\tau _{g}^{2}}})^{^{\prime }}T+\frac{\tau _{g}}{\sqrt{%
k_{n}^{2}+\tau _{g}^{2}}}(k_{g}B+k_{n}N)\right] \\
&&\mp \sinh \theta \left[ (\frac{k_{n}}{\sqrt{k_{n}^{2}+\tau _{g}^{2}}}%
)^{^{\prime }}B+\frac{k_{n}}{\sqrt{k_{n}^{2}+\tau _{g}^{2}}}(-k_{g}T+\tau
_{g}N)\right] +\cosh \theta (k_{n}T+\tau _{g}B).
\end{eqnarray*}%
If we arrange this equality, we obtain%
\begin{gather*}
d^{^{\prime }}=\left( \pm \sinh \theta \left[ (\frac{\tau _{g}}{\sqrt{%
k_{n}^{2}+\tau _{g}^{2}}})^{^{\prime }}+\frac{k_{g}k_{n}}{\sqrt{%
k_{n}^{2}+\tau _{g}^{2}}}\right] +k_{n}\cosh \theta \right) T \\
+\left( \pm \sinh \theta \left[ -(\frac{k_{n}}{\sqrt{k_{n}^{2}+\tau _{g}^{2}}%
})^{^{\prime }}+\frac{k_{g}\tau _{g}}{\sqrt{k_{n}^{2}+\tau _{g}^{2}}}\right]
+\tau _{g}\cosh \theta \right) B.
\end{gather*}%
From Eq (3.4), we have%
\begin{equation*}
\cos h\theta =\mp \sinh \theta \frac{\tau _{g}^{^{\prime
}}k_{n}-k_{n}^{^{\prime }}\tau _{g}+k_{g}(k_{n}^{2}+\tau _{g}^{2})}{%
(k_{n}^{2}+\tau _{g}^{2})^{\tfrac{3}{2}}}.
\end{equation*}%
If the last equality is replaced in the statement of $d^{^{\prime }}$, we get%
\begin{gather*}
d^{^{\prime }}=\pm \sinh \theta \left( 
\begin{array}{c}
\dfrac{\tau _{g}^{^{\prime }}(k_{n}^{2}+\tau _{g}^{2})-\tau
_{g}(k_{n}k_{n}^{^{\prime }}+\tau _{g}\tau _{g}^{^{\prime
}})+k_{g}k_{n}(k_{n}^{2}+\tau _{g}^{2})}{(k_{n}^{2}+\tau _{g}^{2})^{\tfrac{3%
}{2}}} \\ 
+\dfrac{k_{n}k_{n}^{^{\prime }}\tau _{g}-k_{n}^{2}\tau _{g}^{^{\prime
}}-k_{g}k_{n}(k_{n}^{2}+\tau _{g}^{2})}{(k_{n}^{2}+\tau _{g}^{2})^{\tfrac{3}{%
2}}}%
\end{array}%
\right) T \\
\pm \sinh \theta \left( 
\begin{array}{c}
\dfrac{-\kappa _{n}^{^{\prime }}(k_{n}^{2}+\tau
_{g}^{2})+k_{n}(k_{n}k_{n}^{^{\prime }}+\tau _{g}\tau _{g}^{^{\prime
}})+k_{g}\tau _{g}(k_{n}^{2}+\tau _{g}^{2})}{(k_{n}^{2}+\tau _{g}^{2})^{%
\tfrac{3}{2}}} \\ 
+\dfrac{k_{n}^{^{\prime }}\tau _{g}^{2}-k_{n}\tau _{g}\tau _{g}^{^{\prime
}}-k_{g}\tau _{g}(k_{n}^{2}+\tau _{g}^{2})}{(k_{n}^{2}+\tau _{g}^{2})^{%
\tfrac{3}{2}}}%
\end{array}%
\right) B.
\end{gather*}%
By the straight-forward calculation, it follows that $d^{^{\prime }}=0$
namely $d$ is a constant vector.

\begin{theorem}
A unit speed curve $\alpha $ on the spacelike surface $M$ is an isophote
curve with the timelike axis $d$ if and only if%
\begin{equation*}
\psi (s)=\pm \left( \frac{k_{n}^{2}}{(k_{n}^{2}+\tau _{g}^{2})^{\tfrac{3}{2}}%
}(\frac{\tau _{g}}{k_{n}})^{^{\prime }}+\frac{k_{g}}{(k_{n}^{2}+\tau
_{g}^{2})^{\tfrac{1}{2}}}\right) (s)
\end{equation*}%
is a constant function.

\begin{proof}
We can show that $\alpha $ is an isophote curve on the spacelike surface if
and only if the Gaussian mapping along the curve $\alpha $ is a latitude
circle on the Lorentzian unit sphere $S_{1}^{2}$. Hence, if we compute the
Gaussian mapping $N_{\mid _{\alpha }}:I\longrightarrow S_{1}^{2}$ along the
curve $\alpha $, the geodesic curvature of $N_{\mid _{\alpha }}$is $\psi (s)$
as shown below.%
\begin{eqnarray*}
N_{\mid _{\alpha }}^{^{\prime }} &=&k_{n}T+\tau _{g}B \\
N_{\mid _{\alpha }}^{^{\prime \prime }} &=&(k_{n}^{^{\prime }}-k_{g}\tau
_{g})T+(k_{n}k_{g}+\tau _{g}^{^{\prime }})B+(k_{n}^{2}+\tau _{g}^{2})N \\
N_{\mid _{\alpha }}^{^{\prime }}\times N_{\mid _{\alpha }}^{^{\prime \prime
}} &=&\tau _{g}(k_{n}^{2}+\tau _{g}^{2})T+\left( k_{g}(k_{n}^{2}+\tau
_{g}^{2})+k_{n}^{2}(\frac{\tau _{g}}{k_{n}})^{^{\prime }}\right)
N-k_{n}(k_{n}^{2}+\tau _{g}^{2})B,
\end{eqnarray*}%
where $T\times B=N$, $B\times N=T$ and $N\times T=B$. Accordingly, we get%
\begin{eqnarray*}
\kappa &=&\frac{\sqrt{\left \langle N_{\mid _{\alpha }}^{^{\prime }}\times
N_{\mid _{\alpha }}^{^{\prime \prime }},N_{\mid _{\alpha }}^{^{\prime
}}\times N_{\mid _{\alpha }}^{^{\prime \prime }}\right \rangle }}{\left
\Vert N_{\mid _{\alpha }}^{^{\prime }}\right \Vert ^{3}} \\
&=&\frac{\sqrt{\tau _{g}^{2}(k_{n}^{2}+\tau _{g}^{2})^{2}-\left(
k_{g}(k_{n}^{2}+\tau _{g}^{2})+k_{n}^{2}(\dfrac{\tau _{g}}{k_{n}})^{^{\prime
}}\right) ^{2}+k_{n}^{2}(k_{n}^{2}+\tau _{g}^{2})^{2}}}{\sqrt{%
(k_{n}^{2}+\tau _{g}^{2})^{3}}} \\
&=&\sqrt{1-\frac{\left( k_{g}(k_{n}^{2}+\tau _{g}^{2})+k_{n}^{2}(\dfrac{\tau
_{g}}{k_{n}})^{^{\prime }}\right) ^{2}}{(k_{n}^{2}+\tau _{g}^{2})^{3}}}
\end{eqnarray*}%
Let $\overset{-}{k_{g}}$ and $\overset{-}{k_{n}}$ be the geodesic curvature
and the normal curvature of the Gaussian mapping $N_{\mid _{\alpha }}$,
respectively. Since the normal curvature $\overset{-}{k_{n}}=1$ if we
substitute $\overset{-}{k_{n}}$ and $\kappa $ in the following equation, we
obtain the geodesic curvature $\overset{-}{k_{g}}$ as follows.%
\begin{equation*}
\kappa ^{2}=(\overset{-}{k_{n}})^{2}-(\overset{-}{k_{g}})^{2}
\end{equation*}%
\begin{equation*}
\overset{-}{k_{g}}(s)=\psi (s)=\coth \theta =\mp \left( \frac{k_{n}^{2}}{%
(k_{n}^{2}+\tau _{g}^{2})^{\tfrac{3}{2}}}(\frac{\tau _{g}}{k_{n}})^{^{\prime
}}+\frac{k_{g}}{(k_{n}^{2}+\tau _{g}^{2})^{\tfrac{1}{2}}}\right) (s),
\end{equation*}%
where $\theta $ is the angle between the surface normal $N$ and the fixed
vector $d$. Then, the spherical image of isophotes are latitude circles if
and only if $\psi (s)$ is a constant function.
\end{proof}
\end{theorem}

\begin{corollary}
If $\alpha $ is a unit speed isophote curve with the timelike axis $d$ on $M$%
, then $\alpha $ cannot be a silhouette curve.
\end{corollary}

\begin{proof}
Let $\alpha $ be a unit speed isophote curve with the timelike axis $d$ on $%
M $. Then,%
\begin{equation*}
\left \langle N,d\right \rangle =-\cosh \theta \neq 0.
\end{equation*}%
Therefore, by the definition of silhouette curve $\alpha $ cannot be a
silhouette curve.
\end{proof}

\begin{proposition}
If $\alpha $ is a unit speed isophote curve with the timelike axis $d$ on $M$%
, then $\alpha $ cannot be a line of curvature.

\begin{proof}
Let $\alpha $ be a unit speed isophote curve on $M$. By Theorem (1) we have%
\begin{equation*}
\coth \theta =\mp \left( \frac{k_{n}^{2}}{(k_{n}^{2}+\tau _{g}^{2})^{\tfrac{3%
}{2}}}(\frac{\tau _{g}}{k_{n}})^{^{\prime }}+\frac{k_{g}}{(k_{n}^{2}+\tau
_{g}^{2})^{\tfrac{1}{2}}}\right) (s)=constant.
\end{equation*}%
Assume that $\alpha $ is a line of curvature. By applying $\tau _{g}=0$ and
Eq (2.4), we obtain $\coth \theta =\mp \dfrac{k_{g}}{k_{n}}=\mp \dfrac{%
\kappa \sinh \phi }{\kappa \cosh \phi }=\mp \tanh \phi $. The last equation $%
\coth \theta =\mp \tanh \phi $ does not have a solution. For this reason, $%
\alpha $ cannot be a unit speed isophote curve on $M$. This is a
contradiction with respect to our assertion. Hence, $\alpha $ cannot be a
line of curvature.
\end{proof}
\end{proposition}

\begin{corollary}
Let $\alpha $ be a unit speed isophote curve with the timelike axis $d$ on $%
M $. Then, we have the following.\newline
(1) The axis $d$ cannot be perpendicular to the tangent line of $\alpha $.%
\newline
(2) The axis $d$ cannot be perpendicular to the vector $B$ for $\alpha $
with timelike principal normal.
\end{corollary}

\begin{proof}
(1) Suppose that $\alpha $ be a unit speed isophote curve on $M$. By Eq
(3.3), it follows that%
\begin{equation*}
\left \langle T,d\right \rangle =\pm \frac{\tau _{g}}{\sqrt{k_{n}^{2}+\tau
_{g}^{2}}}\sinh \theta .
\end{equation*}%
By the definition of isophote curve, we must have $\sinh \theta \neq 0$ and
also by Proposition (3) we have $\tau _{g}\neq 0.$ Accordingly, it follows
that $\left \langle T,d\right \rangle \neq 0$ in the above equation put
differently the axis $d$ cannot be perpendicular to the tangent line of $%
\alpha $.\newline
(2) Suppose that $\alpha $ be a unit speed isophote curve on $M$. By Eq
(3.3) we have%
\begin{equation*}
\left \langle B,d\right \rangle =\mp \frac{k_{n}}{\sqrt{k_{n}^{2}+\tau
_{g}^{2}}}\sinh \theta .
\end{equation*}%
From Proposition (1) and the definition of isophote curve, we have $%
k_{n}\neq 0$ and $\sinh \theta \neq 0$. From this, we conclude that $%
\left
\langle B,d\right \rangle \neq 0$. Then, the axis $d$ cannot be
perpendicular to the vector $B$.
\end{proof}

\begin{lemma}
Let $\alpha $ be a unit speed spacelike curve in $E_{1}^{3}$ with $\kappa
(s)\neq 0$. Then, $\alpha $ is a slant helix with timelike principal normal
if and only if $\sigma (s)=\left( \dfrac{\kappa ^{2}}{(\kappa ^{2}+\tau
^{2})^{\tfrac{3}{2}}}(\dfrac{\tau }{\kappa })^{^{\prime }}\right) (s)$ is a
constant function $[1]$.
\end{lemma}

\begin{theorem}
Let $\alpha $ be a unit speed isophote curve on $M$. In that case, $\alpha $
is a geodesic if and only if $\alpha $ is a slant helix with the timelike
fixed vector%
\begin{equation*}
d=\pm \dfrac{\tau }{\sqrt{\kappa ^{2}+\tau ^{2}}}\sinh \theta T\pm \dfrac{%
\kappa }{\sqrt{\kappa ^{2}+\tau ^{2}}}\sinh \theta B+\cosh \theta N.
\end{equation*}
\end{theorem}

\begin{proof}
Since $\alpha $ is a geodesic (surface normal $N$ concurs with the principal
normal $n$ along the curve $\alpha $), we have $k_{g}=0$ and thus from Eq
(2.4) $k_{n}=\kappa $ and $\tau _{g}=\tau $. By substituting $k_{g}$ and $%
k_{n}$ in the expression of $\psi (s)$, we follow that%
\begin{equation*}
\psi (s)=\mp \left( \frac{\kappa ^{2}}{(\kappa ^{2}+\tau ^{2})^{\tfrac{3}{2}}%
}(\dfrac{\tau }{\kappa })^{^{\prime }}\right) (s)
\end{equation*}%
is a constant function. By Lemma (1) $\alpha $ is a slant helix. Using Eq
(3.3), the timelike axis of slant helix is obtained as%
\begin{equation*}
d=\pm \dfrac{\tau }{\sqrt{\kappa ^{2}+\tau ^{2}}}\sinh \theta T\pm \dfrac{%
\kappa }{\sqrt{\kappa ^{2}+\tau ^{2}}}\sinh \theta B+\cosh \theta N.
\end{equation*}%
On the contrary, let $\alpha $ be a slant helix with the timelike fixed
vector $d$. From Eq (3.3) we have $k_{n}=\kappa $ and $\tau _{g}=\tau $. So,
the geodesic curvature $k_{g}$ must be zero namely $\alpha $ is a geodesic
on $M$.
\end{proof}

From this time, we will obtain the spacelike fixed vector $d$ for the case
2. If the axis $d$ is spacelike, by the Definition (2) we possess%
\begin{equation}
\left \langle N(u,v),d\right \rangle =\sinh \theta .  \tag{3.5}
\end{equation}%
Just as the case of timelike fixed vector $d$, we get%
\begin{equation*}
\left \langle T,d\right \rangle =-\frac{\tau _{g}}{k_{n}}\left \langle
B,d\right \rangle
\end{equation*}%
If we say $\left \langle B,d\right \rangle =a$, because $\left \langle
d,d\right \rangle =1$ we gather that%
\begin{equation*}
a=\mp \frac{k_{n}}{\sqrt{k_{n}^{2}+\tau _{g}^{2}}}\cosh \theta
\end{equation*}%
Thus the spacelike axis $d$ can be written as%
\begin{equation}
d=\pm \frac{\tau _{g}}{\sqrt{k_{n}^{2}+\tau _{g}^{2}}}\cosh \theta T\mp 
\frac{k_{n}}{\sqrt{k_{n}^{2}+\tau _{g}^{2}}}\cosh \theta B-\sinh \theta N. 
\tag{3.6}
\end{equation}%
Seeing that $\left \langle N^{^{\prime \prime }},d\right \rangle =0$, we
obtain%
\begin{gather}
\coth \theta =\pm \frac{(k_{n}^{2}+\tau _{g}^{2})^{\tfrac{3}{2}}}{%
k_{g}(k_{n}^{2}+\tau _{g}^{2})+(\tau _{g}^{^{\prime }}k_{n}-k_{n}^{^{\prime
}}\tau _{g})}  \tag{3.7} \\
\tanh \theta =\pm \left[ \frac{k_{n}^{2}}{(k_{n}^{2}+\tau _{g}^{2})^{\tfrac{3%
}{2}}}\left( \frac{\tau _{g}}{k_{n}}\right) ^{^{\prime }}+\frac{k_{g}}{%
(k_{n}^{2}+\tau _{g}^{2})^{\tfrac{1}{2}}}\right] .  \notag
\end{gather}%
Applying Eq (3.7) it can be showed that $d^{^{\prime }}=0$ in other words $d$
is a constant vector.

\begin{theorem}
A unit speed spacelike curve $\alpha $ on $M$ is an isophote with the
spacelike axis $d$ if and only if%
\begin{equation*}
\tanh \theta =\omega (s)=\pm \left( \frac{k_{n}^{2}}{(k_{n}^{2}+\tau
_{g}^{2})^{\tfrac{3}{2}}}\left( \frac{\tau _{g}}{k_{n}}\right) ^{^{\prime }}+%
\frac{k_{g}}{(k_{n}^{2}+\tau _{g}^{2})^{\tfrac{1}{2}}}\right) (s)
\end{equation*}%
is a constant function.\newline
The proof of this theorem is the same as Theorem (1).
\end{theorem}

\begin{corollary}
If $\alpha $ is a unit speed isophote curve with the spacelike axis $d$ on $%
M $, then $\alpha $ cannot be a silhouette curve.

\begin{proof}
Let $\alpha $ be a unit speed isophote curve with spacelike axis $d$ on $M$.
Then,%
\begin{equation*}
\left \langle N,d\right \rangle =\sinh \theta .
\end{equation*}%
According to the definition of silhouette curve the surface normal must be
orthogonal to the fixed vector. In the above equality when $\theta =0$, $%
\left \langle N,d\right \rangle =0$. This means that the surface normal $N$
coincides with the axis $d$ but it is not possible to define an isophote
curve like this. Hence, $\alpha $ cannot be a silhouette curve.
\end{proof}
\end{corollary}

\begin{theorem}
Let $\alpha $ be a unit speed isophote curve with the spacelike axis $d$ on $%
M$. Then, $\alpha $ is a geodesic on $M$ if and only if $\alpha $ is a slant
helix with the spacelike fixed vector%
\begin{equation*}
d=\pm \frac{\tau }{\sqrt{\kappa ^{2}+\tau ^{2}}}\cosh \theta T\mp \frac{%
\kappa }{\sqrt{\kappa ^{2}+\tau ^{2}}}\cosh \theta B-\sinh \theta N.
\end{equation*}%
The proof of this theorem like Theorem (2).
\end{theorem}

\begin{proposition}
Let $\alpha $ be a unit speed isophote curve with the spacelike axis $d$ on $%
M$. Then, $\alpha $ is a plane curve provided that $\alpha $ is a line of
curvature on $M$.
\end{proposition}

\begin{proof}
Let $\alpha $ be a unit speed isophote curve on $M$. By Theorem (4) we have%
\begin{equation*}
\tanh \theta =\pm \left( \frac{k_{n}^{2}}{(k_{n}^{2}+\tau _{g}^{2})^{\tfrac{3%
}{2}}}\left( \frac{\tau _{g}}{k_{n}}\right) ^{^{\prime }}+\frac{k_{g}}{%
(k_{n}^{2}+\tau _{g}^{2})^{\tfrac{1}{2}}}\right) (s)=constant.
\end{equation*}%
Assume that $\alpha $ is a line of curvature. By applying $\tau _{g}=0$ and
Eq (2.4), we obtain $\tanh \theta =\pm \dfrac{k_{g}}{k_{n}}=\pm \dfrac{%
\kappa \sinh \phi }{\kappa \cosh \phi }=\pm \tanh \phi $. In consequence, it
concludes that $\phi =\pm \theta $. Since $\phi $ is a constant, using $\tau
_{g}=\tau +\phi ^{^{\prime }}=0$ we get $\tau =0$. Then, $\alpha $ is a
plane curve.
\end{proof}

\begin{theorem}
Let $\alpha $ be a unit speed isophote curve with the spacelike axis $d$ on $%
M$. Then, we have the following.\newline
(1) The axis $d$ is orthogonal to the tangent line of $\alpha $ if and only
if $\alpha $ is a line of curvature on $M$.\newline
(2) The axis $d$ cannot be orthogonal to the vector $B$ for $\alpha $ with
timelike principal normal.
\end{theorem}

\begin{proof}
(1) Let the axis $d$ be orthogonal to the tangent line of $\alpha $. From Eq
(3.6)%
\begin{equation*}
\left \langle T,d\right \rangle =\pm \frac{\tau _{g}}{\sqrt{k_{n}^{2}+\tau
_{g}^{2}}}\cosh \theta =0.
\end{equation*}%
In the equation above, since $\cosh \theta \neq 0$, we get $\tau _{g}=0$.
Consequently, the axis $d$ is orthogonal to the tangent line of $\alpha $ if
and only if $\alpha $ is a line of curvature on $M$.\newline
(2) According to Eq (3.6)%
\begin{equation*}
\left \langle B,d\right \rangle =\mp \frac{k_{n}}{\sqrt{k_{n}^{2}+\tau
_{g}^{2}}}\cosh \theta =0.
\end{equation*}%
From Proposition (1), we have $k_{n}\neq 0$. Also, we have $\cosh \theta
\neq 0$ in the above equation. Then $\left \langle B,d\right \rangle \neq 0$%
, i.e., the axis $d$ cannot be orthogonal to the vector $B$.
\end{proof}

\section{Conclusions}

In this paper, we found the axis (fixed vector) of an isophote curve through
its Darboux frame in $E_{1}^{3}$. Subsequently, we obtained some
characterizations regarding these curves. By using the characterizations, it
is investigated relation between special curves on a spacelike surface and
an isophote. For instance, we viewed the curve which is both isophote and
geodesic on the spacelike surface is a slant helix. Also, we viewed that
isophote curve with the timelike axis cannot also be a line of curvature
while isophote curve with the spacelike axis can be a line of curvature.

\end{document}